\documentclass[12pt,reqno]{amsart}
\usepackage{amsmath, amsfonts, amssymb, amsthm,hyperref,mathtools,array,bm}
\hypersetup{hypertex=true,
	colorlinks=true,
	linkcolor=blue,
	anchorcolor=blue,
	citecolor=blue}
\allowdisplaybreaks[4]
\def\({\bg(}
\def\){\bg)}

\def\Gal{{\rm Gal}}

\def\v{{\bm v}}

\def\alg{{\rm alg}}

\def\Re{{\rm Re}}
\def\B{{\rm B}}
\def\pmod #1{\ ({\rm{mod}}\ #1)}
\def\mod #1{\ {\rm mod}\ #1}
\def\Ack{\medskip\noindent {\bf Acknowledgments}}

\theoremstyle{plain}
\newtheorem{theorem}{Theorem}[section]
\newtheorem{lemma}{Lemma}

\newtheorem{example}{Example}

\theoremstyle{definition}

\theoremstyle{remark}
\newtheorem{remark}{Remark}

\newcommand{\sign}[1]{\mathrm{sign}(#1)}
\makeatletter
\@namedef{subjclassname@2020}{%
	\textup{2020} Mathematics Subject Classification}
\makeatother
\vspace{4mm}

\begin{document}
	\medskip
	\title[On $p$-th cyclotomic field and cyclotomic matrices involving Jacobi sums]
	{On $p$-th cyclotomic field and cyclotomic matrices involving Jacobi sums}
		\author[H.-L. Wu, L.-Y. Wang and H. Pan]{Hai-Liang Wu, Li-Yuan Wang and Hao Pan*}
	
		\address {(Hai-Liang Wu) School of Science, Nanjing University of Posts and Telecommunications, Nanjing 210023, People's Republic of China}
		\email{\tt whl.math@smail.nju.edu.cn}
		
		\address {(Li-Yuan Wang) School of Physical and Mathematical Sciences, Nanjing Tech University, Nanjing 211816, People's Republic of China}
		\email{\tt wly@smail.nju.edu.cn}
		
		\address {(Hao Pan) School of Applied Mathematics, Nanjing University of Finance and Economics, Nanjing 210046, People's Republic of China}
		\email{\tt haopan79@zoho.com}
	
	\keywords{Jacobi sums, Gauss sums, Cyclotomic matrices, Finite field analogues.
		\newline \indent *Corresponding author.}
	
	\begin{abstract}
		Inspired by Weil's classical result on the zeta function of a projective Fermat curve defined over a finite field, in this paper, we investigate some arithmetic properties of the cyclotomic matrix 
		$$\left[J_p(\chi^{ki},\chi^{kj})\right]_{1\le i,j\le n-1},$$
		 where $p\ge3$ is a prime, $1\le k<p-1$ is a divisor of $p-1$ with $p-1=kn$, $\chi$ is a generator of the group of all multiplicative characters of $\mathbb{F}_p$ and $J_p(\chi^{ki},\chi^{kj})$ is the Jacobi sum. For example, let $\zeta_p\in\mathbb{C}$ be a primitive $p$-th root of unity and $P_k(T)$ be the minimal polynomial of the algebraic integer 
		$$\theta_k=\sum_{x\in\mathbb{F}_p,x^k=1}\zeta_p^x$$
		over $\mathbb{Q}$. Then we prove that 
		$$\det \left[J_p(\chi^{ki},\chi^{kj})\right]_{1\le i,j\le n-1}=(-1)^{\frac{(k+1)(n^2-n)}{2}}\cdot n^{n-2}\cdot x_p(k),$$
		where $x_p(k)$ is the coefficient of $T$ in $P_k(T)$. 
	\end{abstract}
	\maketitle

	\section{Introduction}
	\setcounter{lemma}{0}
	\setcounter{theorem}{0}
	\setcounter{equation}{0}
	\setcounter{conjecture}{0}
	\setcounter{remark}{0}
	\setcounter{corollary}{0}
	
	\subsection{Notation}
	
	Throughout this paper, $p$ denotes a prime. Let $\mathbb{F}_p$ be the finite field with $p$ elements, and $\mathbb{F}_p^{\times}$ denotes the multiplicative group of all nonzero elements of $\mathbb{F}_p$. Let $\widehat{\mathbb{F}_p^{\times}}$ be the group of all multiplicative characters of $\mathbb{F}_p$. For any multiplicative character 
	$$A:\ \mathbb{F}_p^{\times}\rightarrow \mathbb{C}^{\times},$$
	we define $A(0)=0$. Also, the inverse character of $A$, denoted by $\bar{A}$, is defined by $\bar{A}(x)=\overline{A(x)}$ for any $x\in\mathbb{F}_p$, where $\bar{z}$ is the complex conjugate of a complex number $z$. As usual, $\varepsilon$ denotes the trivial character. 
	
	Let $\zeta_p=e^{2\pi i/p}\in\mathbb{C}$ be a primitive $p$-th root of unity. Then, for any $A,B\in\widehat{\mathbb{F}_p^{\times}}$ the Gauss sum of $A$ is defined by 
	$$G_p(A):=\sum_{x\in\mathbb{F}_p}A(x)\zeta_p^x,$$
	and the Jacobi sum of $A$ and $B$ is defined by 
	$$J_p(A,B):=\sum_{x\in\mathbb{F}_p}A(x)B(1-x).$$
	
	In addition, for any square matrix $M$ over a field, we use $\det M$ to denote the determinant of $M$, and the entry in the $i$-th row and $j$-th column of $M$ is denoted by $M(i,j)$. For any positive integer $m$, the $m\times m$ identity matrix is denoted by $I_m$.

	\subsection{Background and Motivation} Jacobi sums have extensive applications in number theory, and share many similarities with the arithmetic properties of the beta function. Recall that the beta function is defined by 
	$$\B(x,y)=\int_{0}^{1}t^{x-1}(1-t)^{y-1}dt,$$
	where $x,y\in\mathbb{C}$ with $\Re(x),\Re(y)>0$. Then it is known that (cf. \cite{F,G}) the function
	$$J_p:\ \widehat{\mathbb{F}_p^{\times}}\times\widehat{\mathbb{F}_p^{\times}}\rightarrow \mathbb{Q}(\zeta_{p-1})$$
	defined by $(A,B)\mapsto J_p(A,B)$ is a finite field analogue of the beta function, where $\zeta_{p-1}\in\mathbb{C}$ is a primitive $(p-1)$-th root of unity.
	
	 Normand \cite{Normand} studied many interesting determinants concerning the beta function. For example, by \cite[(4.8)]{Normand} one can easily verify that 
	\begin{equation}\label{Eq. det of Beta function}
		\det\left[\B(i,j)\right]_{1\le i,j\le n}=(-1)^{\frac{n(n-1)}{2}}\prod_{r=0}^{n-1}\frac{(r!)^3}{(n+r)!}.
	\end{equation}
	As a natural finite field analogue of (\ref{Eq. det of Beta function}), it was proved in \cite[Theorem 1.4]{WWP} that 
	\begin{equation}\label{Eq. first finite field analogue of det of Beta function}
		\det \left[J_p(\chi^i,\chi^j)\right]_{1\le i,j\le p-2}=(p-1)^{p-3},
	\end{equation}
	where $p\ge3$ and $\chi$ is a generator of $\widehat{\mathbb{F}_p^{\times}}$. Let  $1\le k<p-1$ be a divisor of $p-1$ with $p-1=kn$. 
	
	As a generalization of (\ref{Eq. first finite field analogue of det of Beta function}), in this paper, we concentrate on the cyclotomic matrices
	\begin{equation}\label{Eq. def. of X}
		X_{p,\chi}(k):=\left[J_p(\chi^{ki},\chi^{kj})\right]_{1\le i,j\le n-1},
	\end{equation}
	and
	\begin{equation}\label{Eq. def. of Y}
		Y_{p,\chi}(k):=\left[J_p(\chi^{ki},\chi^{kj})\right]_{0\le i,j\le n-1}.
	\end{equation}
	Our main task is to determine the explicit values of $\det X_{p,\chi}(k)$ and $\det Y_{p,\chi}(k)$.  To complete this task, we first state a classical result due to Weil. Let 
	$$C_n:=\left\{[x,y,z]\in\mathbb{P}^2(\mathbb{F}_p^{\alg}): x^n+y^n=z^n\right\}$$
	be the projective Fermat curve defined over $\mathbb{F}_p$, where $\mathbb{F}_p^{\alg}$ is an algebraic closure of $\mathbb{F}_p$ and $\mathbb{P}^2(\mathbb{F}_p^{\alg})$ is the projective plane. Consider the zeta function of $C_n$ defined by 
	$$\zeta_{C_n}\left(T\right):=\exp\left(\sum_{m=1}^{+\infty}\frac{N(p^m)\cdot T^m}{m}\right),$$
	where $N(p^m)$ is the number of $\mathbb{F}_{p^m}$-rational points on $C_n$. As $C_n$ is an absolutely irreducible curve of genus $(n-1)(n-2)/2$, by the Weil theorem (cf. \cite[Chapter 18]{gtm84})
	$$\zeta_{C_n}\left(T\right)=\frac{H(T)}{(1-T)(1-pT)},$$
	where 
	\begin{equation}\label{Eq. definition of H(T)}
		H(T)=\prod_{\substack{1\le i,j\le n-1\\ i+j\not\equiv 0\pmod n}}\left(1+J_p(\chi^{ki},\chi^{kj})T\right)\in\mathbb{Z}[T].
	\end{equation}
	
	Motivated by (\ref{Eq. definition of H(T)}), it is reasonable to conjecture that both $\det X_{p,\chi}(k)$ and $\det Y_{p,\chi}(k)$ may be closely related to the coefficients of certain polynomials. In this paper, we will show that both $\det X_{p,\chi}(k)$ and $\det Y_{p,\chi}(k)$ have close relations with the coefficients of the minimal polynomial of certain algebraic integer. In fact, let 
	$$U_k:=\left\{x\in\mathbb{F}_p^{\times}: x^k=1\right\}$$
	be the set of all $k$-th roots of unity over $\mathbb{F}_p$, 
	and recall that $\zeta_p\in\mathbb{C}$ is a primitive $p$-th root of unity. Then we define 
	$$\theta_k:=\sum_{x\in U_k}\zeta_p^x.$$
	As $\mathbb{Q}(\zeta_p)/\mathbb{Q}$ is a cyclic extension, for the divisor $n$ of $p-1$, there is a unique intermediate field $L_n$ of $\mathbb{Q}(\zeta_p)/\mathbb{Q}$ such that $[L_n:\mathbb{Q}]=n$. In the proof of Theorem \ref{Thm. explicit values of det X and det Y}, we will see that $$\mathbb{Q}(\theta_k)=L_n.$$ 
	Noting that $\theta_k$ is an algebraic integer, by the above the minimal polynomial of $\theta_k$ over $\mathbb{Q}$ can be written as 
	\begin{equation}\label{Eq. definition of Pk(T)}
		P_k(T)=T^n+\cdots+x_p(k)T+y_p(k)\in\mathbb{Z}[T].
	\end{equation}
	
	\subsection{Main Theorems} Now we state our main results. 
	
	\begin{theorem}\label{Thm. explicit values of det X and det Y} 
		Let $p\ge3$ be a prime and $1\le k<p-1$ be a divisor of $p-1$ with $p-1=kn$. Then, for any generator $\chi$ of $\widehat{\mathbb{F}_p^{\times}}$, 
		$$\det X_{p,\chi}(k)=(-1)^{\frac{(k+1)(n^2-n)}{2}}\cdot n^{n-2}\cdot x_p(k),$$
		where $x_p(k)$ is the coefficient of $T$ in the minimal polynomial $P_k(T)$. Also, 
		$$\det Y_{p,\chi}(k)=\frac{1}{p}\cdot(-1)^{\frac{(k+1)(n^2-n)}{2}}\cdot n^n\cdot\left(k^2x_p(k)-y_p(k)\right),$$
		where $y_p(k)$ is the constant term of the minimal polynomial $P_k(T)$.
	\end{theorem}
	
	\begin{remark}\label{Rem. Thm. explicit values of det X and det Y} 
		By Theorem \ref{Thm. explicit values of det X and det Y}, we obtain the following identity:
		$$n^2k^2\cdot \det X_{p,\chi}(k)-p\cdot \det Y_{p,\chi}(k)=(-1)^{\frac{(k+1)(n^2-n)}{2}}\cdot n^n\cdot y_p(k).$$
		In the proof of this theorem, we will see that $\det Y_{p,\chi}(k)$ can be decomposed as the sum of two determinants, one involving the coefficient $x_p(k)$ and the other involving the coefficient $y_p(k)$. This makes the combination  $$k^2x_{p}(k)-y_p(k)$$ 
		appear in the formula for $\det Y_{p,\chi}(k)$. 
	\end{remark}
	
	To make our result more explicit, we list some examples here. 
	\begin{example}
		When $p=5$ and $k=2$, by computations we see that 
		$$U_2=\{1,4\}\subseteq\mathbb{F}_5,$$
		and 
		$$\theta_2=\zeta_5+\zeta_5^4=\zeta_5+\overline{\zeta_5}=2\cos\frac{2\pi }{5}=\frac{\sqrt{5}-1}{2}.$$
		Thus, the minimal polynomial of $\theta_2$ is equal to 
		$$P_2(T)=T^2+T-1$$
		with $x_5(2)=1$ and $y_5(2)=-1$. In this case, $X_{5,\chi}(2)=[J_5(\chi^2,\chi^2)]$ is a $1\times 1$ matrix, where $\chi\in\widehat{\mathbb{F}_5^{\times}}$ is a character of order $4$. Hence $\det X_{5,\chi}(2)=J_5(\chi^2,\chi^2)=-1$. 
	\end{example}
	
    \begin{example}
         When $p=7$ and $k=3$, with the help of a computer we have 
    	$$U_3=\{1,2,4\}\subseteq\mathbb{F}_7,$$
    	and
    	$$\theta_3=\zeta_7+\zeta_7^2+\zeta_7^4=\frac{\sqrt{-7}-1}{2}.$$
    	Hence, the minimal polynomial of $\theta_3$ is equal to 
    	$$P_3(T)=T^2+T+2$$
    	with $x_7(3)=1$ and $y_7(3)=2$. In this case, $X_{7,3}=[J_7(\chi^3,\chi^3)]$ is also a $1\times 1$ matrix, where $\chi\in\widehat{\mathbb{F}_7^{\times}}$ is a character of order $6$. Thus, we obtain $\det X_{7,3}=J_7(\chi^3,\chi^3)=1$. 
    \end{example}

	\begin{example}
		When $p=13$ and $k=4$, we have 
		$$U_4=\{1,5,8,12\}\subseteq\mathbb{F}_{13},$$
		and 
		$$\theta_4=\zeta_{13}+\zeta_{13}^5+\zeta_{13}^8+\zeta_{13}^{12}=2\cos\frac{2\pi}{13}+2\cos\frac{10\pi}{13}.$$
		With the help of a computer, we obtain that the minimal polynomial of $\theta_4$ is equal to 
		$$P_4(T)=T^3+T^2-4T+1$$
		with $x_{13}(4)=-4$ and $y_{13}(4)=1$. Since $2$ is a primitive root modulo $13$, we choose a generator $\chi$ of $\widehat{\mathbb{F}_{13}^{\times}}$ with $\chi(2)=e^{2\pi i/12}$. By this we see that 
		$$X_{13,\chi}(4)=\begin{bmatrix}
			J_{13}(\chi^4,\chi^4) & J_{13}(\chi^4,\chi^8)\\
			J_{13}(\chi^8,\chi^4) & J_{13}(\chi^8,\chi^8)
		\end{bmatrix}.$$
		With the help of a computer we obtain $\det X_{13,\chi}(4)=12$. 
	\end{example}

	Inspired by Theorem \ref{Thm. explicit values of det X and det Y}, it is natural to consider the determinants related to Gauss sums. The next theorem demonstrates that these two types of determinants are essentially similar.
	
	\begin{theorem}\label{Thm. tranformation theorem}
		Let $p\ge3$ be a prime and $1\le k<p-1$ be a divisor of $p-1$ with $p-1=kn$. Then, for any generator $\chi$ of $\widehat{\mathbb{F}_p^{\times}}$,
		$$\det\left[G_p(\chi^{ki+kj})\right]_{1\le i,j\le n-1}=(-1)^{\frac{k(n^2-n)}{2}}\cdot \det X_{p,\chi}(k)=(-1)^{\frac{n^2-n}{2}}\cdot n^{n-2}\cdot x_p(k),$$
		where $x_p(k)$ is the coefficient of $T$ in the minimal polynomial $P_k(T)$. In addition,
		\begin{align*}
			\det\left[G_p(\chi^{ki+kj})\right]_{0\le i,j\le n-1}
			&=(-1)^{\frac{k(n^2-n)}{2}}\left(p\det Y_{p,\chi}(k)-(p-1)^2\det X_{p,\chi}(k)\right)\\
			&=(-1)^{\frac{n^2-n+2}{2}}\cdot n^n\cdot y_p(k),
		\end{align*}  
		where $y_p(k)$ is the constant term of the minimal polynomial $P_k(T)$.
	\end{theorem}

	\subsection{Novelty and comparison with previous work} In recent work \cite{WWP}, by using identities related to Jacobi sums and the Cauchy-Binet formula, the authors constructed matrix decompositions of $X_{p,\chi}(1)$ and $X_{p,\chi}(2)$ and obtained the explicit values of $\det X_{p,\chi}(1)$ and $\det X_{p,\chi}(2)$. For example, it was proved that 
	$$\det X_{p,\chi}(2)=\frac{1+(-1)^{(p+1)/2}p}{4}\cdot \left(\frac{p-1}{2}\right)^{(p-5)/2}$$
	for any prime $p\ge5$. However, the method used in \cite{WWP} cannot handle a general divisor $k$ of $p-1$. In addition, for $k\in\{1,2\}$, the exact values of $\det [G_p(\chi^{ki+kj})]_{0\le i,j\le n-1}$ were obtained in 
	\cite{WLWY}. For instance, 
	$$\det\left[G_p(\chi^{2i+2j})\right]_{0\le i,j\le (p-3)/2}=-\left(\frac{p-1}{2}\right)^{(p-1)/2}$$
	whenever $p\ge3$. However, \cite{WWP} does not reveal the connection between such cyclotomic matrices related to Gauss sums and $X_{p,\chi}(k)$ and $Y_{p,\chi}(k)$. Also, the method used in \cite{WLWY} cannot obtain satisfactory results for a general divisor $k$ of $p-1$. 
	
	In this paper, by combining a new method related to almost circulant matrices established in \cite{WW} and the arithmetic properties of the $p$-th cyclotomic field, we will obtain the exact values of $\det X_{p,\chi}(k)$ and $\det Y_{p,\chi}(k)$ for every positive proper divisor $k$ of $p-1$. Moreover, we will reveal the connection between cyclotomic matrices  with Gauss sum entries and $X_{p,\chi}(k)$ and $Y_{p,\chi}(k)$.

	\subsection{Outline of this paper} 	We will prove our main results in Sections 2--3. 
	
	\section{Proof of Theorem \ref{Thm. explicit values of det X and det Y}}
	\setcounter{lemma}{0}
	\setcounter{theorem}{0}
	\setcounter{equation}{0}
	\setcounter{conjecture}{0}
	\setcounter{remark}{0}
	\setcounter{corollary}{0}
	
	Recall that $\varepsilon$ is the trivial multiplicative character of $\mathbb{F}_p$. We begin with the following known results concerning Gauss sums and Jacobi sums (see \cite[Chapters 1--2]{BEK}).
	\begin{lemma}\label{Lem. basic properties on Gauss sums and Jacobi sums}
		Let $p$ be a prime and let $A,B\in\widehat{\mathbb{F}_p^{\times}}$. Then,
		
		{\rm (i)} If $A\neq\varepsilon$ or $B\neq\varepsilon$, then 
		$$J_p(A,B)=\frac{G_p(A)G_p(B)}{\delta_p(AB)G_p(AB)},$$
		where 
		$$\delta_p(AB)=\begin{cases}
			1 & \mbox{if}\ AB\neq\varepsilon,\\
			p & \mbox{if}\ AB=\varepsilon.
		\end{cases}$$
		
		{\rm (ii)} Suppose $A\neq\varepsilon$. Then, 
		\begin{itemize}
			\item $G_p(A)G_p(\bar{A})=pA(-1)$,
			\item $G_p(\bar{A})=A(-1)\overline{G_p(A)}$,
			\item $J_p(A,\bar{A})=-A(-1)$,
			\item $G_p(\varepsilon)=-1$ and $J_p(\varepsilon,\varepsilon)=p-2$. 
		\end{itemize}  
	\end{lemma}
	
	\begin{remark}\label{Rem. Lem. basic properties on Gauss sums and Jacobi sums}
		Recall that the trivial character $\varepsilon$ is defined by 
		$$\varepsilon(x)=\begin{cases}
			1 & \mbox{if}\ x\in\mathbb{F}_p^{\times},\\
			0 & \mbox{otherwise}.
		\end{cases}$$
		Hence, we see that 
		$$G_p(\varepsilon)=\sum_{x\in\mathbb{F}_p}\varepsilon(x)\zeta_p^x=-1+\sum_{x\in\mathbb{F}_p}\zeta_p^x=-1,$$
		and 
		$$J_p(\varepsilon,\varepsilon)=\sum_{x\in\mathbb{F}_p}\varepsilon(x)\varepsilon(1-x)=\sum_{x\in\mathbb{F}_p\setminus\{0,1\}}1=p-2.$$
	\end{remark}

	We now turn to a result concerning some permutations on finite cyclic groups. Let $n\ge2$ be an integer and let $a\in\mathbb{Z}$ with $\gcd(a,n)=1$. Then the map 
	$$x \mod{n}\mapsto ax\mod{n}$$ 
	induces a permutation $\tau_n(a)$ of $\mathbb{Z}/n\mathbb{Z}$. Lerch \cite{Lerch} obtained the following result which determines the sign of $\tau_n(a)$ completely.
	
	\begin{lemma}\label{Lem. Lerch}
		Let notations be as above and let $\sign{\tau_n(a)}$ denote the sign of the permutation $\tau_n(a)$. Then 
		$$\sign{\tau_n(a)}=\begin{cases}
			(\frac{a}{n})   & \mbox{if}\ n\equiv 1\pmod{2},\\
			1                       & \mbox{if}\ n\equiv 2\pmod{4},\\
			(-1)^{(a-1)/2} & \mbox{if}\ n\equiv 0\pmod{4},
		\end{cases}$$
		where $(\frac{\cdot}{n})$ is the Jacobi symbol if $n$ is odd. In particular, 
		$$\sign{\tau_n(-1)}=(-1)^{\frac{(n-1)(n-2)}{2}}.$$
	\end{lemma}

	Before the statement of our next lemma, we briefly introduce some notations for circulant matrices and almost circulant matrices. Let $n\ge2$ be an integer and let 
	$$\v=(a_0,a_1,\cdots,a_{n-1})\in\mathbb{C}^n.$$
	Then the circulant matrix of $\v$ is an $n\times n$ matrix defined by 
	$$C(\v):=[a_{j-i}]_{0\le i,j\le n-1},$$
	where $a_s=a_t$ whenever $s\equiv t\pmod {n}$. More precisely,
	\begin{equation}\label{Eq. definition of circulant matrix}
		C(\v)=[a_{j-i}]_{0\le i,j\le n-1}=\begin{bmatrix}
			a_0         & a_1       & \cdots  & a_{n-2} & a_{n-1}\\
			a_{n-1}   & a_0      & \cdots  & a_{n-3} & a_{n-2}\\
			\vdots    & \vdots  & \ddots  & \vdots    & \vdots\\
			a_2         & a_3      & \cdots  & a_0         & a_1\\
			a_1         & a_2       & \cdots  & a_{n-1} & a_0
		\end{bmatrix}.
	\end{equation}
	On the other hand, the almost circulant matrix $W(\v)$ of $\v$ is an $(n-1)\times (n-1)$ matrix obtained by deleting the first row and the first column of $C(\v)$, that is, 
	\begin{equation}\label{Eq. definition of almost circulant matrix}
		W(\v)=[a_{j-i}]_{1\le i,j\le n-1}=\begin{bmatrix}
			a_0      & \cdots  & a_{n-3} & a_{n-2}\\
			\vdots  & \ddots  & \vdots    & \vdots\\
			a_3      & \cdots  & a_0         & a_1\\
			a_2       & \cdots  & a_{n-1} & a_0
		\end{bmatrix}.
	\end{equation}
	
	The first author and the second author \cite[Theorem 4.1]{WW} proved the following result.
	
	\begin{lemma}\label{Lem. evaluations of almost circulant determinant}
		Let notations be as above and let $\lambda_0,\lambda_1,\cdots,\lambda_{n-1}$ be exactly all the eigenvalues of $C(\v)$. Then 
		$$\det W(\v)=\frac{1}{n}\sum_{l=0}^{n-1}\prod_{r\in[0,n-1]\setminus\{l\}}\lambda_r,$$
		where $[0,n-1]=\{0,1,\cdots,n-1\}$. 
	\end{lemma}
	
	Recall that 
	$$U_k=\{x\in\mathbb{F}_p:\ x^k=1\}$$ and 
	$$\theta_k=\sum_{x\in U_k}\zeta_p^{x}.$$
	We next state the following known result (see \cite[Lemma 3.5 and Lemma 4.8]{WLWY}) concerning the eigenvalues of the matrix $[G_p(\chi^{-ki+kj})^{-1}]_{0\le i,j\le n-1}$.

	\begin{lemma}\label{Lem. eigenvalues of Gp inverse}
		Let $p$ be a prime and $1\le k<p-1$ be a divisor of $p-1$ with $p-1=kn$. For any $b\in\mathbb{F}_p^{\times}/U_k$, let 
		$$\theta_k^{(b)}=\sum_{x\in U_k}\zeta_p^{bx}.$$
		Then for any generator $\chi$ of $\widehat{\mathbb{F}_p^{\times}}$, the complex  numbers 		
		$$\frac{1}{p}-1+\frac{n}{p}\cdot\theta_k^{(b)}\ (\text{where}\  b\in\mathbb{F}_p^{\times}/U_k)$$ 
		are precisely all the eigenvalues of the matrix $[G_p(\chi^{-ki+kj})^{-1}]_{0\le i,j\le n-1}$. 
	\end{lemma}
	
	Now we are in a position to prove our first theorem.
	
	{\noindent{\bf Proof of Theorem \ref{Thm. explicit values of det X and det Y}}.}  (i) Since $\Gal\left(\mathbb{Q}(\zeta_p)/\mathbb{Q}\right)\cong\mathbb{F}_p^{\times}$, we can set 
	$$\Gal\left(\mathbb{Q}(\zeta_p)/\mathbb{Q}\right)=\left\{\sigma_s: s\in\mathbb{F}_p^{\times}\right\},$$
	where $\sigma_s(\zeta_p)=\zeta_p^s$ for any $s\in\mathbb{F}_p^{\times}$.
	
	Since $1,\zeta_p,\cdots,\zeta_p^{p-2}$ is a basis of the extension $\mathbb{Q}(\zeta_p)/\mathbb{Q}$, the numbers  $\zeta_p,\zeta_p^2,\cdots,\zeta_p^{p-1}$ are linearly independent over $\mathbb{Q}$. Thus, for any $\sigma_s\in\Gal\left(\mathbb{Q}(\zeta_p)/\mathbb{Q}\right)$, we have 
	\begin{align*}
		\sigma_s\in\Gal\left(\mathbb{Q}(\zeta_p)/\mathbb{Q}(\theta_k)\right)
		&\Leftrightarrow\sigma_s(\theta_k)=\theta_k\\
		&\Leftrightarrow\sum_{y\in U_k}\zeta_p^{sy}=\sum_{y\in U_k}\zeta_p^y\\
		&\Leftrightarrow s\cdot U_k=\{sy: y\in U_k\}=U_k\\
		&\Leftrightarrow s\in U_k.
	\end{align*}
	This implies that $\Gal\left(\mathbb{Q}(\zeta_p)/\mathbb{Q}(\theta_k)\right)\cong U_k$ and hence  $\Gal\left(\mathbb{Q}(\theta_k)/\mathbb{Q}\right)\cong\mathbb{F}_p^{\times}/U_k$. By this, the minimal polynomial of $\theta_k$ over $\mathbb{Q}$ is equal to
	$$P_k(T)=\prod_{\sigma\in\Gal\left(\mathbb{Q}(\theta_k)/\mathbb{Q}\right)}\left(T-\sigma(\theta_k)\right)=\prod_{b\in\mathbb{F}_p^{\times}/U_k}\left(T-\theta_k^{(b)}\right),$$
	where $\theta_k^{(b)}$ is defined by Lemma \ref{Lem. eigenvalues of Gp inverse}. Therefore, the coefficient of $T$ in $P_k(T)$ is equal to 
	\begin{equation}\label{Eq. xp(k)}
		x_p(k)=(-1)^{n-1}\sum_{b\in\mathbb{F}_p^{\times}/U_k}\prod_{c\neq b}\theta_k^{(c)},
	\end{equation}
	and the constant term of $P_k(T)$ is equal to 
	\begin{equation}\label{Eq. yp(k)}
		y_p(k)=(-1)^n\prod_{b\in\mathbb{F}_p^{\times}/U_k}\theta_k^{(b)}.
	\end{equation}
	
	We now evaluate $\det X_{p,\chi}(k)$. By Lemma \ref{Lem. Lerch} we first obtain 
	\begin{align}\label{Eq. replace j by -j}
		\det X_{p,\chi}(k)
		&=\sign{\tau_n(-1)}\cdot \det \left[J_p(\chi^{-ki},\chi^{kj})\right]_{1\le i,j\le n-1} \notag\\
		&=(-1)^{\frac{(n-1)(n-2)}{2}}\cdot \det \left[J_p(\chi^{-ki},\chi^{kj})\right]_{1\le i,j\le n-1}.
	\end{align}
	By (\ref{Eq. replace j by -j}), it remains to consider $\det[J_p(\chi^{-ki},\chi^{kj})]_{1\le i,j\le n-1}$. Recall that in Lemma \ref{Lem. basic properties on Gauss sums and Jacobi sums} the symbol $\delta_p(\cdot)$ is defined by 
	$$\delta_p(\chi^{-ki+kj})=\begin{cases}
		1 & \mbox{if}\ 1\le i\neq j\le n,\\
		p & \mbox{if}\ 1\le i=j\le n.
	\end{cases}$$
   Applying Lemma \ref{Lem. basic properties on Gauss sums and Jacobi sums}(i), for any $1\le i,j\le n-1$ we have 
   $$J_p(\chi^{-ki},\chi^{kj})=\frac{G_p(\chi^{-ki})\cdot G_p(\chi^{kj})}{\delta_p(\chi^{-ki+kj})\cdot G_p(\chi^{-ki+kj})}.$$
   Hence, 
   \begin{align*}
   	     &\det [J_p(\chi^{-ki},\chi^{kj})]_{1\le i,j\le n-1}\\
   	   =&\det\left[\frac{G_p(\chi^{-ki})\cdot G_p(\chi^{kj})}{\delta_p(\chi^{-ki+kj})\cdot G_p(\chi^{-ki+kj})}\right]_{1\le i,j\le n-1}\\
   	   =&\prod_{1\le i\le n-1}G_p(\chi^{-ki})\cdot \prod_{1\le j\le n-1}G_p(\chi^{kj})\cdot \det\left[\frac{1}{\delta_p(\chi^{-ki+kj})\cdot G_p(\chi^{-ki+kj})}\right]_{1\le i,j\le n-1}\\
   	   =&\prod_{1\le i\le n-1}G_p(\chi^{-ki})G_p(\chi^{ki})\cdot \det\left[\frac{1}{\delta_p(\chi^{-ki+kj})\cdot G_p(\chi^{-ki+kj})}\right]_{1\le i,j\le n-1}.
   \end{align*}
	By Lemma \ref{Lem. basic properties on Gauss sums and Jacobi sums}(ii) we have 
	$$G_p(\chi^{-ki})G_p(\chi^{ki})=p\chi^{ki}(-1)=p(-1)^{ki}$$
	for any $1\le i\le n-1$. Hence,  
	\begin{align}\label{Eq. transform to an almost circulant matrix}
		  \det \left[J_p(\chi^{-ki},\chi^{kj})\right]_{1\le i,j\le n-1}
		&=\prod_{1\le i\le n-1}p\chi^{ki}(-1)\cdot \det\left[\frac{1}{\delta_p(\chi^{-ki+kj})\cdot G_p(\chi^{-ki+kj})}\right]_{1\le i,j\le n-1}\notag\\
		&=(-1)^{\frac{kn(n-1)}{2}}p^{n-1} \det\left[\frac{1}{\delta_p(\chi^{-ki+kj})\cdot G_p(\chi^{-ki+kj})}\right]_{1\le i,j\le n-1}.
	\end{align}
	Let the vector
	\begin{align*}
		\v_p
		&=\left(\frac{1}{\delta_p(\varepsilon) G_p(\varepsilon)},\frac{1}{\delta_p(\chi^k)G_p(\chi^k)},\frac{1}{\delta_p(\chi^{2k})G_p(\chi^{2k})},\cdots,\frac{1}{\delta_p(\chi^{(n-1)k})G_p(\chi^{(n-1)k})}\right)\\
		&=\left(-\frac{1}{p},\frac{1}{G_p(\chi^k)},\frac{1}{G_p(\chi^{2k})},\cdots,\frac{1}{G_p(\chi^{(n-1)k})}\right).
	\end{align*}
	Then, it follows from (\ref{Eq. definition of almost circulant matrix}) that 
	\begin{equation}\label{Eq. almost circulant matrix of vector vp}
		\left[\frac{1}{\delta_p(\chi^{-ki+kj})\cdot G_p(\chi^{-ki+kj})}\right]_{1\le i,j\le n-1}=W(\v_p)
	\end{equation}
	is the almost circulant matrix of $\v_p$. Now by (\ref{Eq. transform to an almost circulant matrix}), it is sufficient to evaluate $\det W(\v_p)$. Let $C(\v_p)$ be the circulant matrix of $\v_p$. Then, by the definition of $\delta_p(\cdot)$ we have 
	\begin{equation}\label{Eq. circulant matrix of vector vp}
		C(\v_p)=\left[\frac{1}{\delta_p(\chi^{-ki+kj})\cdot G_p(\chi^{-ki+kj})}\right]_{0\le i,j\le n-1}=\frac{p-1}{p}I_{n}+\left[\frac{1}{G_p(\chi^{-ki+kj})}\right]_{0\le i,j\le n-1}.
	\end{equation}
	Combining this with Lemma \ref{Lem. eigenvalues of Gp inverse}, we see that the numbers 
	\begin{equation}\label{Eq. eigenvalues of C(vp)}
		\left(-\frac{1}{p}+1\right)+\left(\frac{1}{p}-1+\frac{n}{p}\theta_k^{(b)}\right)=\frac{n}{p}\theta_k^{(b)}\ (\text{where}\ b\in\mathbb{F}_p^{\times}/U_k)
	\end{equation}
	are exactly all the eigenvalues of $C(\v_p)$.  By this and Lemma \ref{Lem. evaluations of almost circulant determinant}, we see that 
	\begin{align*}
		\det W(\v_p)
		&=\det \left[\frac{1}{\delta_p(\chi^{-ki+kj})\cdot G_p(\chi^{-ki+kj})}\right]_{1\le i,j\le n-1}
		\\
		&=\frac{1}{n}\sum_{b\in\mathbb{F}_p^{\times}/U_k}\prod_{c\neq b}\frac{n}{p}\theta_k^{(c)}\\
		&=\frac{n^{n-2}}{p^{n-1}}\sum_{b\in\mathbb{F}_p^{\times}/U_k}\prod_{c\neq b}\theta_k^{(c)}\\
		&=(-1)^{n-1}\cdot \frac{n^{n-2}}{p^{n-1}}\cdot x_p(k),
	\end{align*}
	where the last equality follows from (\ref{Eq. xp(k)}). Combining this with (\ref{Eq. replace j by -j}) and (\ref{Eq. transform to an almost circulant matrix}), we have 
	\begin{align*}
		\det X_{p,\chi}(k)
		&=(-1)^{\frac{(n-1)(n-2)}{2}}\cdot \det \left[J_p(\chi^{-ki},\chi^{kj})\right]_{1\le i,j\le n-1}\\
		&=(-1)^{\frac{(n-1)(n-2)}{2}}\cdot (-1)^{\frac{kn(n-1)}{2}}\cdot p^{n-1}\cdot \det W(\v_p)\\
		&=(-1)^{\frac{(n-1)(n-2)}{2}}\cdot (-1)^{\frac{kn(n-1)}{2}}\cdot p^{n-1}\cdot (-1)^{n-1}\cdot \frac{n^{n-2}}{p^{n-1}}\cdot x_p(k)\\
		&=(-1)^{\frac{(k+1)(n^2-n)}{2}}\cdot n^{n-2}\cdot x_p(k).
	\end{align*}
	
	(ii) We turn to $\det Y_{p,\chi}(k)$. Analogous to (\ref{Eq. replace j by -j}), it follows from Lemma \ref{Lem. Lerch} that 
	\begin{equation}\label{Eq. Y case replace j by -j}
		\det Y_{p,\chi}(k)=(-1)^{\frac{(n-1)(n-2)}{2}}\cdot \det \left[J_p(\chi^{-ki},\chi^{kj})\right]_{0\le i,j\le n-1}.
	\end{equation}
	Thus, we next focus on $\det [J_p(\chi^{-ki},\chi^{kj})]_{0\le i,j\le n-1}$. Observe that 
	\begin{align*}
		\left[J_p(\chi^{-ki},\chi^{kj})\right]_{0\le i,j\le n-1}
		&=\begin{bmatrix}
			p-2                                                 & J_p(\varepsilon,\chi^k)       & \cdots & J_p(\varepsilon,\chi^{(n-1)k})\\
			J_p(\chi^{-k},\varepsilon)         & J_p(\chi^{-k},\chi^k)           & \cdots & J_p(\chi^{-k},\chi^{(n-1)k})\\
			\vdots                                             & \vdots                                     & \ddots & \vdots\\
			J_p(\chi^{-(n-1)k},\varepsilon) & J_p(\chi^{-(n-1)k},\chi^k)  & \cdots & J_p(\chi^{-(n-1)k},\chi^{(n-1)k})
       \end{bmatrix}\\
       &=\begin{bmatrix}
       	   -\frac{1}{p}+(p-2+\frac{1}{p}  )                                               & J_p(\varepsilon,\chi^k)       & \cdots & J_p(\varepsilon,\chi^{(n-1)k})\\
       	   J_p(\chi^{-k},\varepsilon)+0        & J_p(\chi^{-k},\chi^k)           & \cdots & J_p(\chi^{-k},\chi^{(n-1)k})\\
       	   \vdots                                             & \vdots                                     & \ddots & \vdots\\
       	   J_p(\chi^{-(n-1)k},\varepsilon)+0& J_p(\chi^{-(n-1)k},\chi^k)  & \cdots & J_p(\chi^{-(n-1)k},\chi^{(n-1)k})
       \end{bmatrix}.
	\end{align*}
	Hence, $\det [J_p(\chi^{-ki},\chi^{kj})]_{0\le i,j\le n-1}$ can be written as a sum of two determinants. More precisely, 
	\begin{equation}\label{Eq. det Y=det M+det N}
		\det \left[J_p(\chi^{-ki},\chi^{kj})\right]_{0\le i,j\le n-1}=\det M_p+\det N_p,
	\end{equation}
	where $M_p$ and $N_p$ are both $n\times n$ matrices defined by 
	\begin{equation*}
		M_p(i,j)=\begin{cases}
			-1/p                                     & \mbox{if}\ (i,j)=(1,1),\\
			J_p(\chi^{-k(i-1)},\chi^{k(j-1)})  & \mbox{if}\ 1\le i,j\le n\ \text{and}\ (i,j)\neq (1,1),
		\end{cases}
	\end{equation*}
	and 
	\begin{equation*}
		N_p(i,j)=\begin{cases}
			p-2+1/p                                 &  \mbox{if}\ (i,j)=(1,1),\\
			0                                             &  \mbox{if}\ 2\le i\le n\ \text{and}\ j=1,\\
			J_p(\chi^{-k(i-1)},\chi^{k(j-1)})     &  \mbox{otherwise},
		\end{cases}
	\end{equation*}
	respectively. 
	
	We first consider $\det M_p$.  Recall that the circulant matrix 
	$$C(\v_p)=\left[\frac{1}{\delta_p(\chi^{-ki+kj})\cdot G_p(\chi^{-ki+kj})}\right]_{0\le i,j\le n-1}$$
	is defined by (\ref{Eq. circulant matrix of vector vp}). Applying Lemma \ref{Lem. basic properties on Gauss sums and Jacobi sums}, we see that 
	$$M_p(i+1,j+1)=\frac{G_p(\chi^{-ki})G_p(\chi^{kj})}{\delta_p(\chi^{-ki+kj})G_p(\chi^{-ki+kj})}$$
	for any $0\le i,j\le n-1$. Thus, 
	\begin{align*}
		 \det M_p
		&=\det\left[\frac{G_p(\chi^{-ki})G_p(\chi^{kj})}{\delta_p(\chi^{-ki+kj})G_p(\chi^{-ki+kj})}\right]_{0\le i,j\le n-1}\\
		&=\prod_{0\le i\le n-1}G_p(\chi^{-ki})\prod_{0\le j\le n-1}G_p(\chi^{kj})\cdot \det \left[\frac{1}{\delta_p(\chi^{-ki+kj})\cdot G_p(\chi^{-ki+kj})}\right]_{0\le i,j\le n-1}\\
		&=\prod_{0\le i\le n-1}G_p(\chi^{-ki})G_p(\chi^{ki})\cdot \det C(\v_p)\\
		&=\prod_{1\le i\le n-1}p\chi^{ki}(-1)\cdot \det C(\v_p)\\
		&=(-1)^{\frac{kn(n-1)}{2}}\cdot p^{n-1}\cdot \det C(\v_p). 
	\end{align*}
	Combining this with (\ref{Eq. yp(k)}) and (\ref{Eq. eigenvalues of C(vp)}), we obtain 
	\begin{equation}\label{Eq. value of det Mp}
        \det M_p
        =(-1)^{\frac{kn(n-1)}{2}}\cdot p^{n-1}\cdot \prod_{b\in \mathbb{F}_p^{\times}/U_k}\frac{n}{p}\theta_k^{(b)}
        =\frac{1}{p}\cdot (-1)^{\frac{kn(n-1)+2n}{2}}\cdot n^n\cdot y_p(k).
	\end{equation}
	
	We now evaluate $\det N_p$. Note that among the elements in the first column of matrix $N_p$, all entries, except for the element $N_p(1,1)$, are zero. Thus, by (\ref{Eq. replace j by -j}) we obtain 
	\begin{align}\label{Eq. value of det Np}
		\det N_p
		&=\left(\frac{1}{p}+p-2\right)\cdot \det\left[J_p(\chi^{-ki},\chi^{kj})\right]_{1\le i,j\le n-1}\notag\\
		&=\frac{(p-1)^2}{p}\cdot (-1)^{\frac{(n-1)(n-2)}{2}}\cdot \det X_{p,\chi}(k)\notag\\
		&=\frac{(p-1)^2}{p}\cdot (-1)^{\frac{kn(n-1)}{2}+n+1}\cdot n^{n-2}\cdot x_p(k).
	\end{align}
	
	Combining (\ref{Eq. value of det Mp}) and (\ref{Eq. value of det Np}) with (\ref{Eq. det Y=det M+det N}) and (\ref{Eq. Y case replace j by -j}), we finally obtain 
	\begin{align*}
		\det Y_{p,\chi}(k)
		&=(-1)^{\frac{(n-1)(n-2)}{2}}\cdot \det \left[J_p(\chi^{-ki},\chi^{kj})\right]_{0\le i,j\le n-1}\\
		&=(-1)^{\frac{(n-1)(n-2)}{2}}\cdot\left(\det M_p+\det N_p\right)\\
		&=(-1)^{\frac{(n-1)(n-2)}{2}}\cdot (-1)^{\frac{kn(n-1)}{2}+n+1}\cdot n^n\cdot \frac{1}{p}\cdot \left(\frac{(p-1)^2}{n^2}x_p(k)-y_p(k)\right)\\
		&=\frac{1}{p}\cdot(-1)^{\frac{(k+1)(n^2-n)}{2}}\cdot n^n\cdot\left(k^2x_p(k)-y_p(k)\right).
	\end{align*}
	This completes the proof of Theorem \ref{Thm. explicit values of det X and det Y}. \qed 
	
	\section{Proof of Theorem \ref{Thm. tranformation theorem}}
	\setcounter{lemma}{0}
	\setcounter{theorem}{0}
	\setcounter{equation}{0}
	\setcounter{conjecture}{0}
	\setcounter{remark}{0}
	\setcounter{corollary}{0}
	
	Recall that $p\ge3$ is a prime with $p-1=kn$, and $\chi$ is a generator of $\widehat{\mathbb{F}_p^{\times}}$. We begin with the following lemma.
	
	\begin{lemma}\label{Lem. A in the proof of 2nd theorem}
		Let notations be as above. Then,
		\begin{equation}\label{Eq. a in the proof of 2nd Thm.}
			J_p(\chi^{ki},\chi^{kj})=\frac{1}{p}\cdot(-1)^{ki+kj}\cdot G_p(\chi^{ki})\cdot G_p(\chi^{kj})\cdot G_p(\chi^{-ki-kj})
		\end{equation}
		for any $0\le i,j\le n-1$ with $(i,j)\neq (0,0)$. 
	\end{lemma} 
	
	\begin{proof}
		As $0\le i,j\le n-1$ and $(i,j)\neq (0,0)$, one can verify that 
		$$\chi^{ki}\chi^{kj}=\varepsilon\Leftrightarrow i+j=n.$$
		
		Suppose first that $i+j=n$. Then by the above and Lemma \ref{Lem. basic properties on Gauss sums and Jacobi sums}, 
		\begin{equation*}
			\frac{1}{p}(-1)^{p-1}G_p(\chi^{ki})G_p(\chi^{-ki})G_p(\varepsilon)=\frac{1}{p}(-1)^{p}p\chi^{ki}(-1)=-\chi^{ki}(-1)=J_p(\chi^{ki},\chi^{kj}).
		\end{equation*}
		Hence (\ref{Eq. a in the proof of 2nd Thm.}) holds in this case.
		
		Suppose now $i+j\neq n$. Then $\chi^{ki}\chi^{kj}\neq \varepsilon$ in this case. Applying Lemma \ref{Lem. basic properties on Gauss sums and Jacobi sums} again, we obtain 
		\begin{align*}
			J_p(\chi^{ki},\chi^{kj})
			&=\frac{G_p(\chi^{ki})G_p(\chi^{kj})}{G_p(\chi^{ki+kj})}\\
			&=\frac{G_p(\chi^{ki})G_p(\chi^{kj})}{p}\cdot\overline{G_p(\chi^{ki+kj})}\\
			&=\frac{G_p(\chi^{ki})G_p(\chi^{kj})}{p}\cdot\chi^{ki+kj}(-1)\cdot G_p(\chi^{-ki-kj})\\
			&=\frac{1}{p}\cdot(-1)^{ki+kj}\cdot G_p(\chi^{ki})\cdot G_p(\chi^{kj})\cdot G_p(\chi^{-ki-kj}).
		\end{align*}
		Thus, (\ref{Eq. a in the proof of 2nd Thm.}) also holds in this case. This completes the proof.
	\end{proof}
	
	We now prove Theorem \ref{Thm. tranformation theorem}. 
	
	{\noindent{\bf Proof of Theorem \ref{Thm. tranformation theorem}}.} (i) Applying  Lemma \ref{Lem. A in the proof of 2nd theorem} and Lemma \ref{Lem. basic properties on Gauss sums and Jacobi sums}, we have 
	\begin{align*}
		\det X_{p,\chi}(k)
		&=\det\left[\frac{1}{p}(-1)^{ki+kj}G_p(\chi^{ki})G_p(\chi^{kj})G_p(\chi^{-ki-kj})\right]_{1\le i,j\le n-1}\\
		&=\frac{1}{p^{n-1}}\prod_{1\le i\le n-1}(-1)^{ki}G_p(\chi^{ki})\prod_{1\le j\le n-1}(-1)^{kj}G_p(\chi^{kj})\cdot \det\left[G_p(\chi^{-ki-kj})\right]_{1\le i,j\le n-1}\\
		&=\frac{1}{p^{n-1}}\prod_{1\le i\le n-1}(-1)^{ki}G_p(\chi^{ki})\prod_{1\le j\le n-1}(-1)^{-kj}G_p(\chi^{-kj})\cdot \det\left[G_p(\chi^{-ki-kj})\right]_{1\le i,j\le n-1}\\
		&=\frac{1}{p^{n-1}}\prod_{1\le i\le n-1}(-1)^{ki}G_p(\chi^{ki})(-1)^{-ki}G_p(\chi^{-ki})\cdot \det\left[G_p(\chi^{-ki-kj})\right]_{1\le i,j\le n-1}\\
		&=\frac{1}{p^{n-1}}\prod_{1\le i\le n-1}p(-1)^{ki}\cdot \det\left[G_p(\chi^{-ki-kj})\right]_{1\le i,j\le n-1}\\
		&=(-1)^{\frac{k(n^2-n)}{2}}\cdot\det\left[G_p(\chi^{-ki-kj})\right]_{1\le i,j\le n-1}.
	\end{align*}
	Combining this with Lemma \ref{Lem. Lerch}, we obtain 
	\begin{align*}
		\det\left[G_p(\chi^{ki+kj})\right]_{1\le i,j\le n-1}
		&=\sign{\tau_n(-1)}\cdot\sign{\tau_n(-1)}\cdot\det\left[G_p(\chi^{-ki-kj})\right]_{1\le i,j\le n-1}\\
		&=\det\left[G_p(\chi^{-ki-kj})\right]_{1\le i,j\le n-1}\\
		&=(-1)^{\frac{k(n^2-n)}{2}}\cdot \det X_{p,\chi}(k).
	\end{align*}
	By this and Theorem \ref{Thm. explicit values of det X and det Y}, we see that 
	$$\det\left[G_p(\chi^{ki+kj})\right]_{1\le i,j\le n-1}=(-1)^{\frac{k(n^2-n)}{2}}\cdot \det X_{p,\chi}(k)=(-1)^{\frac{n^2-n}{2}}\cdot n^{n-2}\cdot x_p(k).$$
	
	(ii) Note that 
	\begin{align*}
		Y_{p,\chi}(k)
		&=\begin{bmatrix}
			p-2                                                 & J_p(\varepsilon,\chi^k)       & \cdots & J_p(\varepsilon,\chi^{(n-1)k})\\
			J_p(\chi^{k},\varepsilon)         & J_p(\chi^{k},\chi^k)           & \cdots & J_p(\chi^{k},\chi^{(n-1)k})\\
			\vdots                                             & \vdots                                     & \ddots & \vdots\\
			J_p(\chi^{(n-1)k},\varepsilon) & J_p(\chi^{(n-1)k},\chi^k)  & \cdots & J_p(\chi^{(n-1)k},\chi^{(n-1)k})
		\end{bmatrix}\\
		&=\begin{bmatrix}
			-\frac{1}{p}+(p-2+\frac{1}{p}  )                                               & J_p(\varepsilon,\chi^k)       & \cdots & J_p(\varepsilon,\chi^{(n-1)k})\\
			J_p(\chi^{k},\varepsilon)+0        & J_p(\chi^{k},\chi^k)           & \cdots & J_p(\chi^{k},\chi^{(n-1)k})\\
			\vdots                                             & \vdots                                     & \ddots & \vdots\\
			J_p(\chi^{(n-1)k},\varepsilon)+0& J_p(\chi^{(n-1)k},\chi^k)  & \cdots & J_p(\chi^{(n-1)k},\chi^{(n-1)k})
		\end{bmatrix}.
	\end{align*}
	Analogous to (\ref{Eq. det Y=det M+det N}), we have the decomposition
	\begin{equation}\label{Eq. det Y=det R+det S in 2nd theorem}
		\det Y_{p,\chi}(k)=\det R_p+\det S_p,
	\end{equation}
	where $R_p$ and $S_p$ are both $n\times n$ matrices defined by 
	\begin{equation*}
		R_p(i,j)=\begin{cases}
			-1/p                                                 & \mbox{if}\ (i,j)=(1,1),\\
			J_p(\chi^{k(i-1)},\chi^{k(j-1)})  & \mbox{if}\ 1\le i,j\le n\ \text{and}\ (i,j)\neq (1,1),
		\end{cases}
	\end{equation*}
	and 
	\begin{equation*}
		S_p(i,j)=\begin{cases}
			p-2+1/p                                             &  \mbox{if}\ (i,j)=(1,1),\\
			0                                                         &  \mbox{if}\ 2\le i\le n\ \text{and}\ j=1,\\
			J_p(\chi^{k(i-1)},\chi^{k(j-1)})     &  \mbox{otherwise},
		\end{cases}
	\end{equation*}
	respectively. 
	
	We first consider $\det R_p$. Observe that, when $(i,j)=(0,0)$, the right hand side of the equality (\ref{Eq. a in the proof of 2nd Thm.}) is $-1/p$, which coincides with $R_p(1,1)$. Thus, applying Lemma \ref{Lem. A in the proof of 2nd theorem}, we obtain 
	$$R_p(i+1,j+1)=\frac{1}{p}\cdot(-1)^{ki+kj}\cdot G_p(\chi^{ki})\cdot G_p(\chi^{kj})\cdot G_p(\chi^{-ki-kj})$$
	for any $0\le i,j\le n-1$. This implies that 
	\begin{align*}
		\det R_p
		&=\det\left[\frac{1}{p}(-1)^{ki+kj}G_p(\chi^{ki})G_p(\chi^{kj})G_p(\chi^{-ki-kj})\right]_{0\le i,j\le n-1}\\
		&=\frac{1}{p^{n}}\prod_{0\le i\le n-1}(-1)^{ki}G_p(\chi^{ki})\prod_{0\le j\le n-1}(-1)^{kj}G_p(\chi^{kj})\cdot \det\left[G_p(\chi^{-ki-kj})\right]_{0\le i,j\le n-1}\\
		&=\frac{1}{p^{n}}\prod_{0\le i\le n-1}(-1)^{ki}G_p(\chi^{ki})\prod_{0\le j\le n-1}(-1)^{-kj}G_p(\chi^{-kj})\cdot \det\left[G_p(\chi^{-ki-kj})\right]_{0\le i,j\le n-1}\\
		&=\frac{1}{p^{n}}\prod_{0\le i\le n-1}(-1)^{ki}G_p(\chi^{ki})(-1)^{-ki}G_p(\chi^{-ki})\cdot \det\left[G_p(\chi^{-ki-kj})\right]_{0\le i,j\le n-1}\\
		&=\frac{1}{p^{n}}\prod_{1\le i\le n-1}p(-1)^{ki}\cdot \det\left[G_p(\chi^{-ki-kj})\right]_{0\le i,j\le n-1}\\
		&=\frac{1}{p}\cdot(-1)^{\frac{k(n^2-n)}{2}}\cdot\det\left[G_p(\chi^{-ki-kj})\right]_{0\le i,j\le n-1}.
	\end{align*}
	By this and Lemma \ref{Lem. Lerch}, we immediately obtain 
	\begin{align}\label{Eq. det R and det G in 2nd theorem}
		\det\left[G_p(\chi^{ki+kj})\right]_{0\le i,j\le n-1}
		&=\sign{\tau_n(-1)}^2\cdot \det\left[G_p(\chi^{-ki-kj})\right]_{0\le i,j\le n-1}\notag\\
		&=(-1)^{\frac{k(n^2-n)}{2}}\cdot p \cdot \det R_p.
	\end{align}
	
	On the other hand, by the definition of matrix $S_p$ we see that 
	\begin{equation}\label{Eq. det S and det X in 2nd theorem}
		\det S_p=\left(p-2+\frac{1}{p}\right)\cdot \det X_{p,\chi}(k)=\frac{(p-1)^2}{p}\cdot \det X_{p,\chi}(k).
	\end{equation}
	
	Now combining (\ref{Eq. det R and det G in 2nd theorem}) and (\ref{Eq. det S and det X in 2nd theorem}) with (\ref{Eq. det Y=det R+det S in 2nd theorem}) and using Theorem \ref{Thm. explicit values of det X and det Y}, we finally obtain 
	\begin{align*}
		\det\left[G_p(\chi^{ki+kj})\right]_{0\le i,j\le n-1}
		&=(-1)^{\frac{k(n^2-n)}{2}}\cdot p \cdot \det R_p\\
		&=(-1)^{\frac{k(n^2-n)}{2}}\cdot p \cdot \left(\det Y_{p,\chi}(k)-\det S_p\right)\\
		&=(-1)^{\frac{k(n^2-n)}{2}}\cdot p \cdot \left(\det Y_{p,\chi}(k)-\frac{(p-1)^2}{p}\det X_{p,\chi}(k)\right)\\
		&=(-1)^{\frac{k(n^2-n)}{2}}\cdot \left(p\det Y_{p,\chi}(k)-(p-1)^2\det X_{p,\chi}(k)\right)\\
		&=(-1)^{\frac{n^2-n+2}{2}}\cdot n^n\cdot y_p(k).
	\end{align*}
	This completes the proof of Theorem \ref{Thm. tranformation theorem}. \qed 
	
	\Ack\ The authors would like to thank the referees for careful reading and helpful comments. These valuable suggestions have significantly improved the quality of this paper.
	
	This work was supported by the National Natural Science Foundation of China (Grant Nos. 12671009, 12201291 and 12471312). The first author was also supported by the Natural Science Foundation of the Higher Education Institutions of Jiangsu Province (Grant No. 25KJB110010). The third author was also supported by the 333 High-level Talents Project of Jiangsu Province (2022-3-16-471).

\end{document}